\documentclass[letterpaper, 10 pt, conference]{ieeeconf}  % Comment this line out
\IEEEoverridecommandlockouts                              % This command is only
\usepackage[latin1]{inputenc}
\usepackage{amsmath}
\usepackage{amssymb}
\usepackage{epsfig}
\usepackage{psfig}
\usepackage{psfrag}
\usepackage{graphics}
\usepackage{color}
\usepackage{subfigure}

\newcommand{\Ebf}{{\mathbf E}}

\newcommand{\Pbf}{{\mathbf P}}

\newcommand{\Xbf}{{\mathbf X}}

\newcommand{\xbf}{{\mathbf x}}

\newcommand{\zbf}{{\mathbf z}}

\newcommand{\Ccal}{{\mathcal {C}}}

\newcommand{\Gcal}{{\mathcal {G}}}
\newcommand{\Ical}{{\mathcal {I}}}
\newcommand{\Pcal}{{\mathcal {P}}}
\newcommand{\Ucal}{{\mathcal {U}}}
\newcommand{\Vcal}{{\mathcal {V}}}
\newcommand{\vcal}{{\mathcal {V}}}
\newcommand{\Xcal}{{\mathcal {X}}}

\newcommand{\rb}{{\mathbb R}}

\newcommand{\beq}{\begin{equation}}
\newcommand{\eeq}{\end{equation}}

\newcommand{\bnul}{\begin{enumerate}[a)]}
\newcommand{\enul}{\end{enumerate}}

\newcommand{\op}[1]{\operatorname{#1}}

\newcommand{\fin}{\hspace*{\fill}~$\blacksquare$}

\allowdisplaybreaks

\title{\LARGE \bf
Optimal input design for non-linear dynamic systems: a graph theory approach
}

\author{Patricio E. Valenzuela, Cristian R. Rojas and Håkan Hjalmarsson % <-this % stops a space
\thanks{This work was supported in part by the Swedish Research Council under contracts 621-2011-5890 and 621-2009-4017, and in part by the European Research Council under the advanced grant LEARN, contract 267381.}% <-this % stops a space
\thanks{Automatic Control Lab and  ACCESS Linnaeus Center, School of Electrical Engineering,
KTH--Royal Institute of Technology, SE-100 44 Stockholm, Sweden.
(e-mail: {\tt \small \{pva, crro, hjalmars\}@kth.se}.)
}%
}

\begin{document}

\maketitle
\thispagestyle{empty}
\pagestyle{empty}

%\addtolength{\textheight}{-3cm}   % This command serves to balance the column lengths
%                                  % on the last page of the document manually. It shortens
%                                  % the textheight of the last page by a suitable amount.
%                                  % This command does not take effect until the next page
%                                  % so it should come on the page before the last. Make
%                                  % sure that you do not shorten the textheight too much.

%%%%%%%%%%%%%%%%%%%%%%%%%%%%%%%%%%%%%%%%%%%%%%%%%%%%%%%%%%%%%%%%%%%%%%%%%%%%%%%%
\begin{abstract}
In this article a new algorithm for the design of stationary input sequences for system identification is presented. The stationary input signal is generated by optimizing an approximation of a scalar function of the information matrix, %. The approximation is
 based on stationary input sequences generated from prime cycles, which describe the set of finite Markov chains of a given order. This method can be used for solving input design problems for nonlinear systems. In particular it can handle amplitude constraints on the input. %, and when amplitude constraints are imposed on the input signal.  %The feasible set is constrained to input sequences adopting a finite set of values.
 Numerical examples show that the new algorithm is computationally attractive and that is consistent with previously reported results. % previously reported in the literature.
\end{abstract}

%%%%%%%%%%%%%%%%%%%%%%%%%%%%%%%%%%%%%%%%%%%%%%%%%%%%%%%%%%%%%%%%%%%%%%%%%%%%%%%%
\begin{keywords}
System identification, input design, Markov chains.
\end{keywords}

%%%%%%%%%%%%%%%%%%%%%%%%%%%%%%%%%%%%%%%%%%%%%%%%%%%%%%%%%%%%%%%%%%%%%%%%%%%%%%%%
\section{INTRODUCTION}\label{sec: 1}
Input design concerns the generation of an input signal to maximize the information obtained from an experiment. %optimize a cost function.  The cost function depends on the application where the input signal will be employed. One
 Some of the first contributions in this area have been introduced in \cite{Cox1958,goopay76}. %An extension of input design for dynamical systems is presented in GOODWIN-PAYNE.
 From the roots of \cite{Cox1958,goopay76}, many contributions on the subject have been developed (see \cite{fedorov1972,Whittle1973,hildebrand2003,gevers2005} and the references therein).%for a survey on input design for dynamical systems).

In the case of dynamic systems, the input design problem can be formulated as the optimization of a cost function related to the model to be identified. The results in this area are mainly focused on linear systems. In \cite{jansson2005,lindqvist2000}, linear matrix inequalities (LMI) are employed to solve the input design problem. % when model specifications are needed to include in the problem. A Markov chain approach is presented
 In \cite{brighenti2009}, the input signal is modelled as the output of a Markov chain. %The input sequence is designed by optimizing the transition probabilities of the Markov chain.
  Robust input design is covered in \cite{Rojas2007}, where the input signal is designed to optimize a cost function over the feasible set of the true parameters. %Extensions of input design to closed-loop have been proposed in \cite{gevers2011,jansson2005b,hildebrand2009}. A connection of input design problem with closed-loop control is introduced in \cite{hjalmarsson2005}.
  A time domain approach for input design for system identification is developed in \cite{Suzuki2007}. The results presented above (except for \cite{brighenti2009}) design an input signal without amplitude constraints. %no constraints on its amplitude.
   However, in practical applications, amplitude constraints on the input signal are required due to physical and/or performance limitations. Therefore, input design with amplitude constraints needs more analysis. %requires investigation. %it is necessary to develop a method suitable to include

In recent years, we see growing %an
 interest to extend the results of input design to nonlinear systems. An approach to input design for nonlinear systems by using the knowledge of linear systems is presented in \cite{hjalmarsson2007}. Input design for structured nonlinear systems has been introduced in \cite{vincent2009}. An approach of input design for a particular class of nonlinear systems is presented in \cite{larsson2010}. A particle filter method for input design for nonlinear systems is presented in \cite{gopaluni2011}. The results presented allow to design input signals when the system contains nonlinear functions, but the constraints on the system dynamics and the computational cost required to solve the problem are the main limitations of %to use
  these results. Therefore, it is necessary to develop a method for input design suited for a wide class of nonlinear models %input design including a wide class of nonlinear systems
  and requiring low computational effort. % to solve the problem.

In this article we present a method to solve input design problems with %by considering
 amplitude limitations. The proposed technique also includes nonlinear systems with more general structures than those presented in \cite{larsson2010}. The method designs an input signal which %by optimizing a scalar cost function of the information matrix. The resulting input signal
  is restricted to %can adopt
   a finite set of values, and it is a realization of the optimal stationary process. Since the problem is solved over the set of stationary processes, the feasible set needs to be described by %using
  basis functions. However, finding the basis functions for this set is a hard task. This drawback is solved by using ideas from graph theory \cite{zaman1983,johnson1975,tarjan1972}. By deriving the prime cycles of the de Brujin's graph associated to the feasible set, we can express any element in the set as a convex sum of a finite number of elements. The information matrices associated to these elements can be approximated by a simple average, which reduces the computational costs compared to the method in \cite{gopaluni2011}. A nice feature of this approach is that, if the cost function is convex, the optimization problem can be solved by using convex tools, even if the system is nonlinear. The numerical examples show that this method is consistent with the results presented in \cite{larsson2010}, and that it can be successfully applied to solve input design problems with amplitude limitations.

%One might argue that the method developed here is not suitable for input design, since the true parameters of the system are needed. However, this method can be implemented in a recursive approach, where the estimated parameters in the $r$-th iteration are used to design the input sequence for the $(r+1)$-th estimation. Recursive algorithms are beyond %of
% the scope of this paper, and it will be addressed in a future work.

As with most optimal input design methods, the one proposed in this contribution relies on knowledge of the true system. This difficulty can be overcome by implementing a robust experiment design scheme on top of it \cite{Rojas2007} or via an adaptive procedure, where the input signal is re-designed as more information is being collected from the system \cite{rojas2011adaptive}. Due to space limitations, however, we will not address these issues in the present paper.

The rest of this paper is organized as follows. Section \ref{sec: 1a} presents basic concepts in graph theory. Section \ref{sec: 2} introduces the input design problem. Section \ref{sec: 3} presents the newly proposed method to compute an optimal input signal. In Section \ref{sec: 4} some numerical examples are presented. Finally, Section \ref{sec: 5} concludes the paper. %presents conclusions. % and ideas for future work in the subject.

\subsubsection*{\textbf{Remark}}
In the sequel, we denote by $\rb$ the real set, by $\rb ^{p}$ the set of real $p$-dimensional vectors, and by $\rb^{r\times s}$ the set of real $r\times s$ matrices. The expected value and the probability measure are denoted by $\Ebf \{ \cdot \}$, and $\Pbf \{ \cdot \}$, respectively. Finally, $\det$ and $\op{tr}$ stand for the determinant and the trace functions, respectively.
%%%%%%%%%%%%%%%%%%%%%%%%%%%%%%%%%%%%%%%%%%%%%%%%%%%%%%%%%%%%%%%%%%%%%%%%%%%%%%%%
\section{PRELIMINARIES ON GRAPH THEORY}\label{sec: 1a}
%The results we will introduce in the coming sections require the knowledge of basic concepts in graph theory.
The purpose in this section is to provide a brief background on graph theory to understand the discussion in the next sections. The definitions presented here come from \cite[pp. 77]{johnson1975}.

A \emph{directed graph} $\Gcal _{\Vcal} := (\Vcal ,\Xcal)$ consists of a nonempty and finite set of vertices (or nodes) $\Vcal$ and a set $\Xcal$ of ordered pairs of distinct vertices called \emph{edges}. A \emph{path} in $\Gcal _\Vcal$ is a sequence of vertices $p_{vu} := (v=v_1,\, v_2 ,\,  \ldots , \, v_k = u)$ such that $(v_i ,\, v_{i+1}) \in \Xcal$ for $i \in \{1, \, \ldots , \, k-1 \}$. A \emph{cycle} is a path in which the first and last vertices are identical. %A path is \emph{elementary} if no vertex appears twice.
A cycle is elementary if no vertex but the first and last appears twice. Two elementary cycles are distinct if one is not a cyclic permutation of the other.

%%%%%%%%%%%%%%%%%%%%%%%%%%%%%%%%%%%%%%%%%%%%%%%%%%%%%%%%%%%%%%%%%%%%%%%%%%%%%%%%
\section{PROBLEM FORMULATION}\label{sec: 2}
\begin{figure}[thpb]
      \begin{center}
    \begin{picture}(0,0)%
    \includegraphics{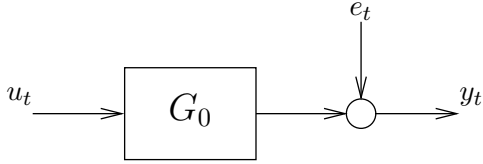}%
    \end{picture}%
    %
    %  Created by WinFIG version 4.9
    %  METADATA <version>1.0</version>
    %
    \setlength{\unitlength}{2171sp}%
    \begingroup\makeatletter\ifx\SetFigFont\undefined%
    \gdef\SetFigFont#1#2#3#4#5{%
      \reset@font\fontsize{#1}{#2pt}%
      \fontfamily{#3}\fontseries{#4}\fontshape{#5}%
      \selectfont}%
    \fi\endgroup%
    \begin{picture}(4919,1955)(2379,-2783)
    %  METADATA <id>7</id>
    \put(6151,-1111){\makebox(0,0)[b]{\smash{{\SetFigFont{11}{13.2}{\rmdefault}{\mddefault}{\updefault}{\color[rgb]{0,0,0}$e_t$}%
    }}}}
    %  METADATA <id>8</id>
    \put(4201,-2311){\makebox(0,0)[b]{\smash{{\SetFigFont{14}{16.8}{\rmdefault}{\mddefault}{\updefault}{\color[rgb]{0,0,0}$G_0$}%
    }}}}
    %  METADATA <id>9</id>
    \put(2401,-2086){\makebox(0,0)[rb]{\smash{{\SetFigFont{11}{13.2}{\rmdefault}{\mddefault}{\updefault}{\color[rgb]{0,0,0}$u_t$}%
    }}}}
    %  METADATA <id>10</id>
    \put(7276,-2086){\makebox(0,0)[lb]{\smash{{\SetFigFont{11}{13.2}{\rmdefault}{\mddefault}{\updefault}{\color[rgb]{0,0,0}$y_t$}%
    }}}}
    \end{picture}%
      \end{center}
      \caption{Block diagram of a dynamic system (possibly nonlinear).}
      \label{fig: 1}
   \end{figure}

Consider the single-input, single-output time invariant system depicted in Figure \ref{fig: 1}. Here, $G_0$ is a dynamic system (possibly nonlinear), $\{e_t\}$ is a white noise sequence with zero mean and variance $\lambda _e$, $u_t \in \rb$ is the input and $y_t \in \rb$ is the measured output. We will assume that we have a model structure for $G_0$. Notice that we assume that the noise $e_t$ enters only at the output. % is known.

The objective in this article is to design an input signal $\Ucal _{n_{\op{seq}}} := (u_{n_{\op{seq}}}, \, \ldots , \, u_1 )$ as a realization of a stationary process, such that the system
\begin{align}
\label{eq: prob1}
y_t &= G_0(\Ucal _{t}^\infty) +e_t\, , \\
\label{eq: prob1a}
\Ucal _{t}^\infty &:= (u_{t}, \, \ldots , \, u_{-\infty} ) \, ,
\end{align}
can be identified  with maximum accuracy as defined by %by optimizing a functional cost. The most common choice for the functional cost is
 a scalar function of the Fisher information matrix $\Ical _F$ \cite{ljung99}. %, which will be used in this paper. If we assume that the prediction error method is used to identify \eqref{eq: prob1}, then
  $\Ical _F$ can be computed as %depends on the identification method, the model structure, the input sequence $U_{n_{\op{seq}}}$ and on the noise sequence $E_{n_{\op{seq}}} := \{e_{(n_{\op{seq}}-1)}, \, \ldots , \, e_0 \}$. In this article, the identification method is the prediction error. Therefore, the model structure is given by
\beq
\label{eq: prob3}
\Ical _F := \dfrac{1}{\lambda _e} \Ebf \left\{ \sum_{t=1} ^{n_{\op{seq}}}\psi(\theta _0) \psi( \theta_0) ^T\right\} \, ,
\eeq
where
\begin{align}
\label{eq: prob2a}
\psi(\theta _0) &:= \left. \dfrac{d \, \hat{y}(t)}{d \theta} \right| _{\theta = \theta _0}\, , \\
\label{eq: prob2b}
\hat{y}(t) &:= G(\Ucal _t;\, \theta) \, ,
\end{align}
and $\theta , \, \theta _0 \in \Theta \subset \rb ^m$. The expected value in \eqref{eq: prob3} is with respect to the realizations of $\Ucal _{n_{\op{seq}}}$. %and $\Ecal _{n_{\op{seq}}}:= \{e_{(n_{\op{seq}}-1)}, \, \ldots , \, e_0 \}$.
In addition, the result introduced in \eqref{eq: prob3}-\eqref{eq: prob2b} assumes that there exists a $\theta _0 \in \Theta$ such that $G(\Ucal _t; \, \theta _0) = G_0(\Ucal _t)$ \cite{ljung99}, i.e., that there is no undermodelling; we will make this assumption in the sequel.
%\beq
%\label{eq: prob2}
%\hat{y}(t) = G(U_t;\, \theta) \, ,
%\eeq
%Equation \eqref{eq: prob2} is the one-step predictor of \eqref{eq: prob1}.
%by optimizing a cost function dependent on $U_{n_{\op{seq}}}$ and the identification method. The most common choice for the cost function is a scalar function of

Equation \eqref{eq: prob2b} does not depend on the noise realization. Therefore, we can rewrite \eqref{eq: prob3} as
%\beq
%\label{eq: prob4}
%\Ical _F = \dfrac{1}{\lambda _e} \int _{\Ucal _t \in \rb ^t} \psi(\theta _0) \psi( \theta_0) ^T\, p(\Ucal _t) \, d \Ucal _t\, ,
%\eeq
\beq
\label{eq: prob4}
\Ical _F = \dfrac{1}{\lambda _e} \int _{\Ucal _{n_{\op{seq}}} \in \rb ^{n_{\op{seq}}}} \sum_{t=1} ^{n_{\op{seq}}} \psi(\theta _0) \psi( \theta_0) ^T\, d P(\Ucal _{n_{\op{seq}}})\, ,
\eeq
where $P(\Ucal _{n_{\op{seq}}})$ is the cumulative distribution function of $\Ucal _{n_{\op{seq}}}$. %Since we are interested in designing a stationary input sequence of length $n_{\op{seq}}$, we have in \eqref{eq: prob4} that $t=n_{\op{seq}}-1$.

We note that \eqref{eq: prob4} depends on $P(\Ucal _{n_{\op{seq}}})$. Therefore, %to obtain an optimal stationary input sequence $\Ucal _{n_{\op{seq}}} ^{\op{opt}}$,
 the input design problem is to find a cumulative distribution function $P^{\op{opt}}(\Ucal _{n_{\op{seq}}})$ which optimizes a scalar function of \eqref{eq: prob4}. We define this scalar function as $h : \, \rb ^{m \times m} \rightarrow \, \rb$. To obtain the desired results, $h$ must be a matrix nondecreasing function \cite[pp. 108]{boyvan04}. Different choices of $h$ have been proposed in the literature \cite{Rojas2007}. Some examples for $h$ are $h = \det$, and $h = -\op{tr} \{(\cdot)^{-1}\}$. In this work, we leave to the user the selection of $h$.

Since $P^{\op{opt}}(\Ucal _{n_{\op{seq}}})$ has to be a stationary cumulative distribution function, the optimization must be constrained to the set
\begin{multline}
\label{eq: prob5}
\Pcal := \left\{ F:\, \rb ^{n_{\op{seq}}} \rightarrow \rb | \, F(\xbf) \geq 0 , \, \forall \xbf \in \rb ^{n_{\op{seq}}}; \,  \right. \\
\left. F \text{ is monotone non-decreasing}\right. ; \, \\
\left. \lim _ {\substack{x_i \rightarrow \infty \\ i = \{1, \, \ldots , \, n_{\op{seq}}\}\\ \xbf = (x_1 , \, \ldots , \, x_{n_{\op{seq}}})}} F(\xbf) = 1 ;\,\right. \\
\left. \int _{v \in \rb} dF(v, \, \zbf) = \int _{v \in \rb} dF(\zbf , \, v) \, , \forall \zbf \in \rb ^{(n_{\op{seq}}-1)} \right\} \, .
\end{multline}
The last condition in \eqref{eq: prob5} (with slight abuse of notation) guarantees that $F \in \Pcal$ is the cumulative distribution function of a stationary sequence \cite{zaman1983}.

To simplify our analysis, we will assume that $u_t$ can only adopt a finite number $c_{\op{seq}}$ of values. We define this set of values as $\Ccal$. With the previous assumption, we can define the following subset of $\Pcal$:% as
\begin{multline}
\label{eq: prob6}
\Pcal _\Ccal := \left\{ f:\, \Ccal ^{n_{\op{seq}}} \rightarrow \rb | \, f(\xbf) \geq 0 , \, \forall \xbf \in \Ccal ^{n_{\op{seq}}}; \,  \right. \\
\left. \sum _ {\xbf \in \Ccal ^{n_{\op{seq}}}} f(\xbf) = 1; \right.  \\
\left. \sum _{v \in \Ccal} f(v, \, \zbf) = \sum _{v \in \Ccal} f(\zbf , \, v) \, , \forall \zbf \in \Ccal ^{(n_{\op{seq}}-1)} \right\} \, .
\end{multline}
The set introduced in \eqref{eq: prob6} will be used to constrain the probability mass function $p(\Ucal _{n_{\op{seq}}})$.

The discussion presented in this section can be summarized as
\subsubsection*{\textbf{Problem 1}}
Design an optimal input signal $\Ucal _{n_{seq}} ^{\op{opt}} \in \Ccal ^{n_{\op{seq}}}$ as a realization from $p^{\op{opt}}(\Ucal _{n_{seq}})$, where% such that
\beq
\label{eq: prob7}
p^{\op{opt}}(\Ucal _{n_{seq}}) := \arg \max _{p \in \Pcal _\Ccal} h(\Ical _F(p)) \, ,
\eeq
where $h : \, \rb ^{m \times m} \rightarrow \, \rb$ is a matrix nondecreasing function,
\beq
\label{eq: prob8}
\Ical _F(p) = \dfrac{1}{\lambda _e} \sum _{\Ucal _{n_{\op{seq}}} \in \Ccal ^{n_{\op{seq}}}}  \sum_{t=1} ^{n_{\op{seq}}} \psi(\theta _0) \psi( \theta_0) ^T\, p(\Ucal _{n_{\op{seq}}}) \, ,
\eeq
and $\psi(\theta _0) \in \rb ^m$ is defined as in \eqref{eq: prob2a}-\eqref{eq: prob2b}.
\fin

A solution for this problem will be discussed in the next section.
%%%%%%%%%%%%%%%%%%%%%%%%%%%%%%%%%%%%%%%%%%%%%%%%%%%%%%%%%%%%%%%%%%%%%%%%%%%%%%%%
\section{INPUT DESIGN VIA GRAPH THEORY}\label{sec: 3}
%To get the optimal input sequence, we need to be able to solve Problem 1. However,
Problem 1 is hard to solve explicitly. The main issues are %drawbacks for Problem 1 are
\begin{enumerate}
\item We need to describe the elements in the set $\Pcal _\Ccal$ as a linear combination of basis functions, and %is infinite dimensional, and
\item the sum in \eqref{eq: prob8} is of dimension $n_{\op{seq}}$, where $n_{\op{seq}}$ could be potentially very large.
\end{enumerate}
%The previous elements
These issues make Problem 1 computationally intractable. Therefore, we need to develop an approach to solve this problem by using a computational feasible method.

Since $n_{\op{seq}}$ could be large, Problem 1 can be unfeasible to solve. To address this, we restrict the memory of the stationary process $u_t$, i.e., we consider only finite stationary sequences of length, say, $n_m$. %In this section, you may then change "n_SEQ" to "n".

To address the first issue, we notice that $\Pcal _\Ccal$ is a convex set. In particular, $\Pcal _\Ccal$ is a polyhedron \cite[pp. 31]{boyvan04}. %Since $\Pcal _\Ccal$ is a polyhedron, then
 Hence, any element of $\Pcal _\Ccal$ can be described as a convex combination of the extreme points of $\Pcal _\Ccal$ \cite[pp. 24]{boyvan04}. Therefore, if we define $\Vcal_{\Pcal _\Ccal}$ as the set of all the extreme points of $\Pcal _\Ccal$, composed by $n_\vcal$ elements, then for all $f \in \Pcal _\Ccal$ we have
\beq
\label{eq: prob9}
f = \sum _{i=1}^{n_{\vcal}} \alpha _i \, v_i \, ,
\eeq
where $\alpha _i \geq 0$, $i \in \{1, \ldots , \, n_\vcal\}$, %and
\beq
\label{eq: prob10}
\sum _{i=1}^{n_{\vcal}} \alpha _i = 1 \, ,
\eeq
and $v_i  \in \Vcal_{\Pcal _\Ccal}$, for all $i \in \{1, \ldots , \, n_\vcal\}$.

Equation \eqref{eq: prob9} says that all the elements in $\Pcal _\Ccal$ can be described by using $n_{\vcal}$ elements in the set $\Vcal_{\Pcal _\Ccal}$. %Therefore, if we are able to find all the elements in the set $\Vcal_{\Pcal _\Ccal}$, then we can overcome drawback 1), which reduces significantly the computational cost of solving Problem 1.

To find all the elements in $\Vcal_{\Pcal _\Ccal}$, we need to shift our focus to graph theory. Indeed, we can analyze the set $\Ccal ^{n_{m}}$ as follows. $\Ccal ^{n_{m}}$ is composed of $(c_{\op{seq}} )^{n_{m}}$ elements. Each element in $\Ccal ^{n_{m}}$ can be viewed as one node in a graph. In addition, the transitions among the elements in $\Ccal ^{n_{m}}$ are given by the feasible values of $u_{t+k+1}$ when we move from $(u_{t+k} , \ldots , \, u_k)$ to $(u_{t+k+1} , \ldots , \, u_{k+1})$, for all integers $k \geq 0$. %Since the transitions among the elements in $\Ccal ^{n_{\op{seq}}}$ are clearly defined,
 The edges among the elements in $\Ccal ^{n_{m}}$ denote the possible transitions between the states, represented by the nodes of the graph. %must be consistent with this assumption.
  Figure \ref{fig: 2} illustrates this idea, when $c_{\op{seq}}= 2$, $n_{m}=2$, and $\Ccal = \{0,\, 1\}$. From this figure we can see that, if we are in node $(0,\, 1)$ at time $t$, then we can only end at node $(1, \, 0)$ or $(1, \, 1)$ at time $t+1$. %If the last one is not satisfied, then the graph is not derived from $\Ccal ^{n_{\op{seq}}}$. %a set with sequences.
   \begin{figure}[t]%hpb]
      \centering
      \input{graph.pstex_t}
      \caption{Example of graph derived from $\Ccal ^{n_{m}}$, with $c_{\op{seq}}=2$, $n_{m}=2$, and $\Ccal := \{0, \, 1\}$.} %(we omit self-edges in $(0,\, 0)$ and $(1, \, 1)$ to improve legibility).}
      \label{fig: 2}
   \end{figure}

The idea to use graph theory to find all the elements in $\Vcal_{\Pcal _\Ccal}$ is related with the concept of prime cycles. In graph theory, a \emph{prime cycle} is an elementary cycle whose set of nodes do not have a proper subset which is an elementary cycle \cite[pp. 678]{zaman1983}. It has been proved that the prime cycles of a stationary graph can describe all the elements in the set $\Vcal_{\Pcal _\Ccal}$ \cite[Theorem 6]{zaman1983}. In other words, each prime cycle defines one element $v_i \in \Vcal_{\Pcal _\Ccal}$. Furthermore, each $v_i$ corresponds to a uniform distribution whose support is the set of elements of its prime cycle, for all $i \in \{1, \ldots , \, n_\vcal\}$ \cite[pp. 681]{zaman1983}. Therefore, the elements in $\Vcal_{\Pcal _\Ccal}$ can be described by finding all the prime cycles associated to the stationary graph $\Gcal _{\Ccal ^{n_{m}}}$ drawn from $\Ccal ^{n_{m}}$.

It is known that all the prime cycles associated to $\Gcal _{\Ccal ^{n_{m}}}$ can be derived from the elementary cycles associated to $\Gcal _{\Ccal ^{(n_{m}-1)}}$ \cite[Lemma 4]{zaman1983}. %An \emph{elementary circuit} is a path where no vertex but the first and last appears twice \cite[pp. 77]{johnson1975}.
 In the literature there are many algorithms for finding all the elementary cycles in a graph. For the examples in Section \ref{sec: 4}, we have used the algorithm presented in \cite[pp. 79--80]{johnson1975} complemented with the one proposed in \cite[pp. 157]{tarjan1972}.

\begin{figure}[t]%hpb]
      \centering
      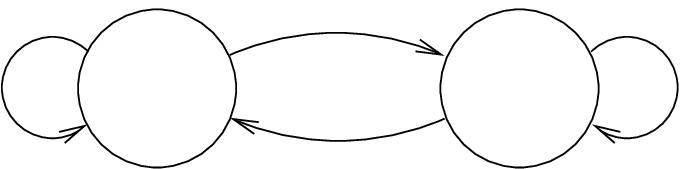
      \caption{Example of graph derived from $\Ccal ^{n_{m}}$, with $c_{\op{seq}}=2$, $n_{m}=1$, and $\Ccal := \{0, \, 1\}$.}
      \label{fig: 2a}
   \end{figure}

Once all the elementary cycles of $\Gcal _{\Ccal ^{(n_{m}-1)}}$ are found, we can find all the prime cycles associated to $\Gcal _{\Ccal ^{n_{m}}}$ by using the idea introduced in \cite[Lemma 4]{zaman1983}. To illustrate this, we consider the graph depicted in Figure \ref{fig: 2a}. One elementary cycle for this graph is given by $(0,\, 1, \, 0)$. %Therefore
 Using Lemma 4 in \cite{zaman1983}, the elements of one prime cycle for the graph $\Gcal _{\Ccal ^{2}}$ are obtained
 as a concatenation of the elements in the elementary cycle $(0,\, 1, \, 0)$. Hence, the prime cycle in $\Gcal _{\Ccal ^{2}}$ associated to this elementary cycle is given by $((0,\, 1), \, (1,\, 0), \, (0,\, 1))$.

With all the prime cycles clearly defined for $\Gcal _{\Ccal ^{n_{m}}}$, then all the elements in the set $\Vcal_{\Pcal_\Ccal}$ are found. Hence, %Since we know all the elements in $\Vcal_{\Pcal_\Ccal}$,
we can use \eqref{eq: prob9} to describe all the elements in $\Pcal _\Ccal$. Thus, the solution described here presents a computationally feasible method to address the first issue.

%With the method presented to solve drawback 1), we can also solve drawback 2).
Since we know the distribution $v_i$ for each prime cycle, with $v_i  \in \Vcal _{\Pcal _\Ccal}$, we can generate an input signal $\{u_t ^i\}_{t=0} ^{t=N}$ drawn from $v_i $, so that
\begin{align}
\nonumber
\Ical _F ^{(i)} &:= \dfrac{1}{\lambda _e} \sum _{\Ucal _{n_{m}} \in \Ccal ^{n_{m}}} \sum_{t=1}^{n_m} \psi(\theta _0) \psi( \theta_0) ^T\, v_i(\Ucal _{n_{m}}) \\
\label{eq: prob11}
&\approx \dfrac{1}{\lambda _e  \, N} \sum _{t=1} ^{N} \psi(\theta _0) \psi( \theta_0) ^T \, ,
\end{align}
for all $i \in \{1, \ldots , \, n_\Vcal\}$, and $N$ sufficiently large\footnote{Note that $N$ is the number of Monte Carlo simulations to compute \eqref{eq: prob11}, and it is not necessarily equal to the length of the experiment $n_{\op{seq}}$.} (in relation to the length of the prime cycles). Notice that $\psi(\theta _0)$ depends implicitly on $\{u_t ^i\}_{t=0} ^{t=N}$ through \eqref{eq: prob2a}-\eqref{eq: prob2b}. Furthermore, each $\Ical _F ^{(i)}$ is associated to the $i$-th prime cycle, for all $i \in \{1, \ldots , \, n_\Vcal\}$.

As an example of how to generate $\{u_t ^i\}_{t=0} ^{t=N}$ from $v_i$, we use the graph depicted in Figure \ref{fig: 2}. One prime cycle for this graph is given by $((0,\, 1), \, (1,\, 0), \, (0, \, 1))$. Therefore, the sequence $\{u_t ^i\}_{t=0} ^{t=N}$ is given by taking the last element of each node, i.e., $\{u_t ^i\}_{t=0} ^{t=N} = \{1,\, 0 , \, 1, \, 0 , \, \ldots, \,  ((-1)^N+1)/2\}$.

The approximation of each $\Ical _F ^{(i)}$ given by \eqref{eq: prob11} reduces the sum \eqref{eq: prob8} from dimension $n_{\op{seq}}$ to dimension 1. This simplification reduces significantly the computation effort to obtain \eqref{eq: prob8}. With this approach, issue 2) is also addressed.

To summarize, the proposed method for input design of signals in $\Ccal ^{n_{m}}$ can be described as follows
\begin{enumerate}
\item Compute all the elementary cycles of $\Gcal _{\Ccal ^{(n_{m}-1)}}$ by using, e.g., \cite[pp. 79--80]{johnson1975}, \cite[pp. 157]{tarjan1972}.
\item Compute all the prime cycles of $\Gcal _{\Ccal ^{n_{m}}}$ from the elementary cycles of $\Gcal _{\Ccal ^{(n_{m}-1)}}$ as explained above (c.f. \cite[Lemma 4]{zaman1983}).
\item Generate the input signals $\{u_t ^i\}_{t=0} ^{t=N}$ from the prime cycles of $\Gcal _{\Ccal ^{n_{m}}}$, for each $i \in \{1, \, \ldots, \, n_\Vcal\}$.
\item For each $i \in \{1, \, \ldots, \, n_\Vcal\}$, approximate $\Ical _F ^{(i)}$ by using \eqref{eq: prob11}.
\item Define $\gamma := \{\alpha _1 , \ldots , \, \alpha _{n_\Vcal}\} \in \rb ^{n_{\Vcal}}$. Find $\gamma ^{\op{opt}} := \{\alpha _1 ^{\op{opt}}, \ldots , \, \alpha _{n_\Vcal} ^{\op{opt}}\}$ by solving an approximation of Problem 1, given by
\beq
\label{eq: prob12}
\gamma ^{\op{opt}} := \arg \max _{\gamma \in \rb ^{n_{\Vcal}}} h(\Ical _F ^{\op{app}}(\gamma)) \, ,
\eeq
where
\begin{align}
\label{eq: prob12a}
\Ical _F ^{\op{app}}(\gamma) &:=  \sum _{i=1} ^{n_{\Vcal}} \alpha _i \, \Ical _F ^{(i)} \, , \\
\sum _{i=1} ^{n_{\Vcal}} \alpha _i &= 1 \, , \\
\label{eq: prob13}
\alpha _i &\geq 0 \, , \text{ for all } i \in \{1, \ldots , \, n_{\Vcal}\} \, ,
\end{align}
and $\Ical _F ^{(i)}$ is given by \eqref{eq: prob11}, for all $i \in \{1, \ldots , \, n_{\Vcal}\}$.
\end{enumerate}
The procedure mentioned above computes $\gamma ^{\op{opt}}$ to describe the optimal probability density function $p^{\op{opt}}(\Ucal _{n_{m}})$ using \eqref{eq: prob9}. Notice that $\Ical _F ^{\op{app}}(\gamma)$ in \eqref{eq: prob12a} is linear in the decision variables. Therefore, for a suitable choice of $h$, the problem \eqref{eq: prob12}-\eqref{eq: prob13} becomes convex.

On the other hand, notice that the steps (1)-(3) mentioned above are independent of the system for which the experiment is designed. Therefore, once steps (1)-(3) are computed, then can be reused to design input sequences for different systems.

 To obtain an input signal from $p^{\op{opt}}(\Ucal _{n_{m}})$, we need to compute a Markov chain associated to the elements in $\Ccal ^{n_{m}}$. %For that reason, we need to compute the stationary probability of each sequence in $\Ccal ^{n_{\op{seq}}}$. Once we computed the stationary probabilities,
  We can find one transition matrix $A \in \rb ^{(c_{\op{seq}})^{n_{m}}\times (c_{\op{seq}})^{n_{m}}}$ for the equivalent Markov chain
\beq
\label{eq: prob14}
\Pi _{k+1} := A \, \Pi _k \, ,
\eeq
by using algorithms presented in the literature (e.g., Metropolis-Hastings algorithm \cite{hastings1970monte,boyd2004}). Notice that each entry of $\Pi _k \in \rb ^{(c_{\op{seq}})^{n_{m}}}$ in \eqref{eq: prob14} represents one element in $\Ccal ^{n_{m}}$. To use the algorithms presented in \cite{hastings1970monte,boyd2004} we need to determine the stationary probabilities of each element in $\Ccal ^{n_{m}}$, which can be computed as follows. %. The stationary probabilities of each element in $\Ccal ^{n_{\op{seq}}}$ can be computed as follows.
 We know that each vertex in $\Vcal _{\Pcal_{\Ccal}}$ has a uniform distribution with support equal to the set of input vectors in the associated prime cycle. Therefore, the stationary probability of each $\xbf \in \Ccal ^{n_{m}}$ is given by
\beq
\label{eq: prob15}
\Pbf \{ \Xbf = \xbf \} = \sum _{i=1} ^{n_{\Vcal}} \alpha _i ^{\op{opt}} \, v_i (\xbf) \, .
\eeq
Equation \eqref{eq: prob15} can be used to construct $\Pi ^{\op{s}} \in \rb ^{(c_{\op{seq}})^{n_{m}}}$, where each entry in $\Pi ^{\op{s}}$ is associated to the stationary probability of one element in $\Ccal ^{n_{m}}$. Given $\Pi ^{\op{s}}$, we can %use \cite{boyd2004} to
 find one matrix $A$ such that %concludes the method proposed to solve Problem 1.
\beq
\label{eq: prob14a}
\Pi ^{\op{s}} = A\, \Pi ^{\op{s}} \, .
\eeq
Finally, the transition matrix $A$ can be used to compute the input sequence by running the Markov chain with a random initial state $\Pi _0$.
%%%%%%%%%%%%%%%%%%%%%%%%%%%%%%%%%%%%%%%%%%%%%%%%%%%%%%%%%%%%%%%%%%%%%%%%%%%%%%%%
\section{NUMERICAL EXAMPLES}\label{sec: 4}
The previous section described a method to compute a solution for Problem 1. In this section we will show that the method is consistent with reported algorithms in the literature.
   \begin{figure}[t]%hpb]
      \centering
      \includegraphics[width = 0.43\textwidth]{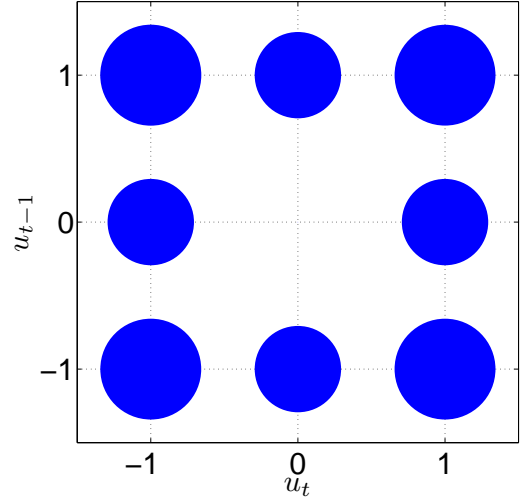}
      \caption{Plot with the stationary probabilities for the optimal input signal of Example 1. The radius of each disc is proportional to the probability of the state $\{u_t, \, u_{t-1}\}$.}
      \label{fig: 3}
   \end{figure}
\subsubsection*{\textbf{Example 1}} In this example we will solve the input design problem for the system in Figure \ref{fig: 1}, with
\beq
\label{eq: prob16}
G_0(\Ucal _t) = G_1(q,\theta) \, u_t + G_2(q,\theta) \, u_t ^2 \, ,
\eeq
    \begin{figure*}[!t]
      \centering
      \subfigure[]{
                    \includegraphics[width = 0.45\textwidth]{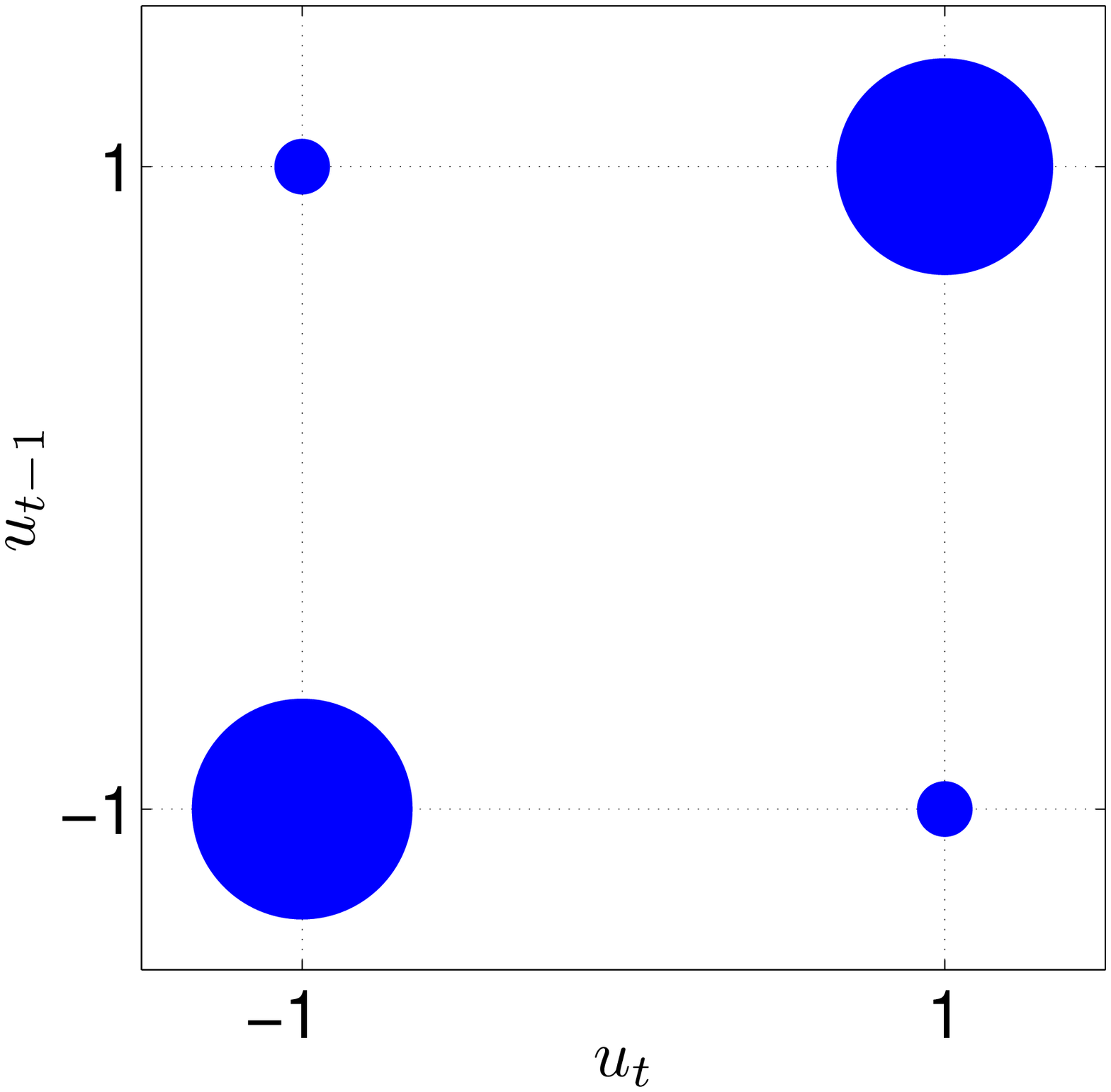}
                    \label{fig: 41}
      }\qquad
      \subfigure[]{
                    \includegraphics[width = 0.45\textwidth]{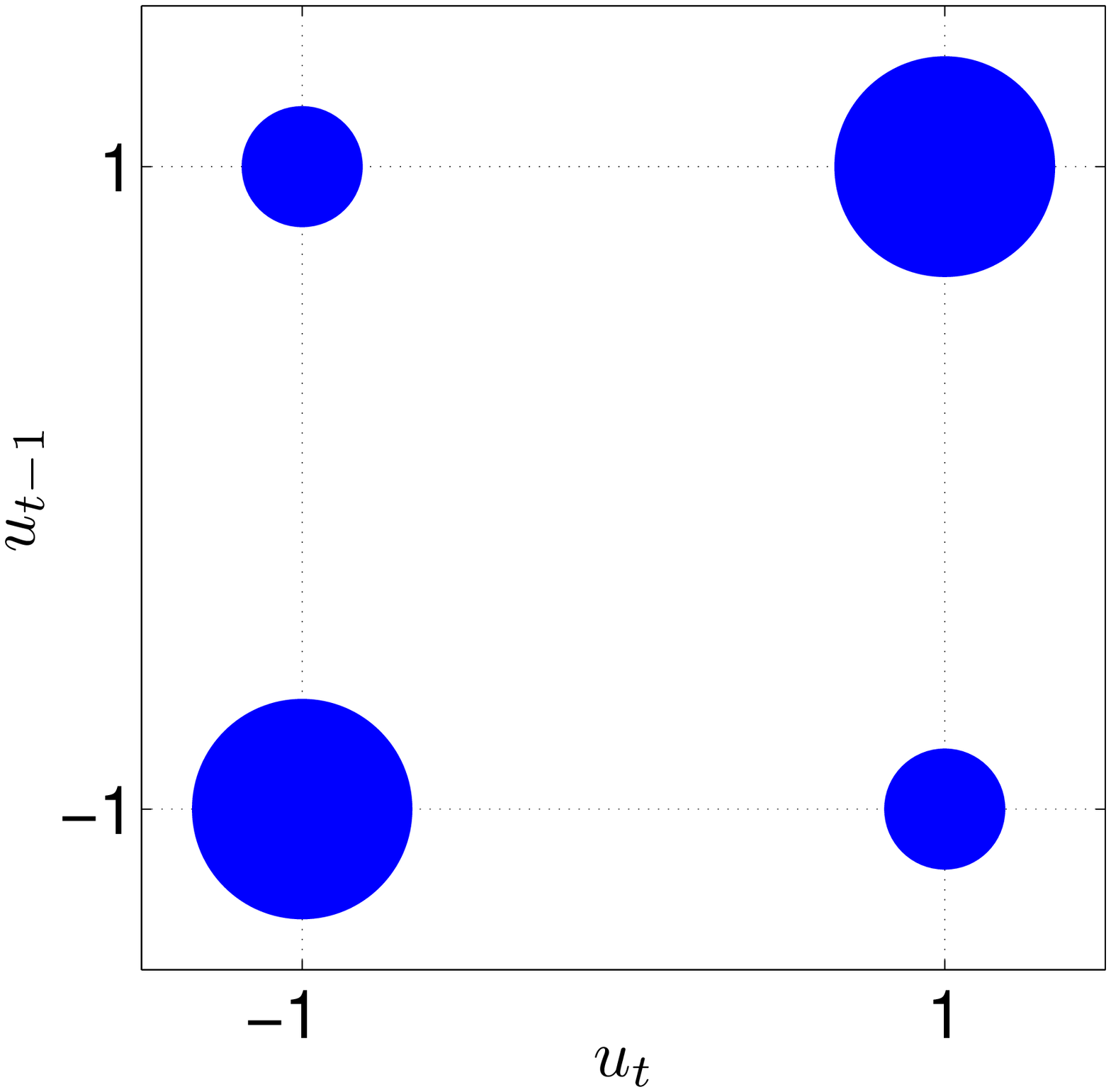}
                    \label{fig: 42}
      }
      \caption{Plot with the stationary probabilities for the optimal input signal in Example 2. The radius of each disc is proportional to the probability of the state $\{u_t, \, u_{t-1}\}$. Figure \ref{fig: 41}: $h(\cdot) = -\op{tr} \{(\cdot)^{-1}\}$. Figure \ref{fig: 42}: $h(\cdot) = \det (\cdot)$.}
      \label{fig: 4}
   \end{figure*}
where
\begin{align}
\label{eq: prob17}
G_1(q,\theta) &= \theta _1 +\theta _2 \, q^{-1} \, , \\
\label{eq: prob18}
G_2(q,\theta) &= \theta _3 +\theta _4 \, q^{-1}  \, ,
\end{align}
and $q$ denotes the shift operator, i.e., $q ^{-1}\, u_t := u_{t-1}$. We assume that $e_t$ is Gaussian white noise with variance $\lambda _e = 1$. This system has been introduced as an example in \cite{larsson2010}.

We will solve Problem 1 by considering $h(\cdot) = \det (\cdot)$, and a ternary sequence ($c_{\op{seq}}=3$) of length $n_{m}=2$. For this example, we define $\Ccal := \{-1, \, 0 , \, 1\}$.

To solve \eqref{eq: prob12}-\eqref{eq: prob13} we consider $N= 5\cdot 10^3$ in \eqref{eq: prob11}. The implementation of \eqref{eq: prob12}-\eqref{eq: prob13} was made in \verb+Matlab+ by using \verb+cvx+ toolbox \cite{boyd_cvxuser2013}.

The simulation results give an optimal cost $\det(\Ical _F ^{\op{app}}) = 0.1796$ (c.f. $\det(P^{-1})=0.18$ for the same example in \cite{larsson2010}). Figure \ref{fig: 3} shows the optimal stationary probabilities for each state $\{u_t, \, u_{t-1}\}$ (c.f. Figure 4(a) in \cite{larsson2010})\footnote{The use of disc plots to represent the optimal input in Figure~\ref{fig: 3} is considered to ease comparison with the results in \cite{larsson2010}, where this visual representation is used.}. The results presented here show that the proposed method is consistent with previous results in the literature \cite{larsson2010}, when $G_0$ is defined as \eqref{eq: prob16}-\eqref{eq: prob18}.
\fin

Example 1 shows that this method is equivalent to the method introduced in \cite{larsson2010} when $G_0$ has a nonlinear FIR-type structure.

The results in this article can be also employed when amplitude constraints are considered in the input sequence by forcing $u_t$ to belong to a finite alphabet. The next example shows an application in that direction. %. as we can see in the next example.
\subsubsection*{\textbf{Example 2}}  In this example we consider the mass-spring-damper system introduced in \cite{brighenti2009}. The continuous input $u$ is the force applied to the mass and the output $y$ is the mass position. The continuous-time system is described by the transfer function
\beq
\label{eq: prob19}
G_0(s) = \dfrac{\frac{1}{m}}{s^2 +\frac{c}{m}\, s +\frac{k}{m}} \, ,
\eeq
with $m = 100$ [Kg], $k = 10$ [N/m], and $c = 6.3246$ [Ns/m]. This choice results in the natural frequency $\omega _n = 0.3162$ [rad/s], and the damping $\xi = 0.1$. The noise $e_t$ is white with zero mean and variance $\lambda _e = 10^{-4}$. The system \eqref{eq: prob19} is sampled by using a zero-order-hold with sampling period $T_s = 1$ [s]. This gives the discrete-time system % which results in
\beq
\label{eq: prob20}
G_0(\Ucal _t) = \dfrac{4.86\cdot 10^{-3}\, q^{-1} +4.75\cdot 10^{-3} q^{-2}}{1 -1.84\, q^{-1} +0.94 q^{-2}}  \, u_t\, .
\eeq
As a model, we define
\beq
\label{eq: prob21}
G(\Ucal _t; \, \theta) = \dfrac{\theta _1 \, q^{-1} +\theta_2 q^{-2}}{1 +\theta_3\, q^{-1} +\theta _4 q^{-2}}  \, u_t\, ,
\eeq
where
\beq
\label{eq: prob22}
\theta = \begin{bmatrix} \theta_1 & \theta_2 & \theta _3 & \theta _4 \end{bmatrix}^T\, .
\eeq
We will solve Problem 1 for two cost functions: $h(\cdot) = -\op{tr} \{(\cdot)^{-1}\}$ and $h(\cdot) = \det (\cdot)$, subject to a binary sequence ($c_{\op{seq}}=2$) of length $n_{m}=2$. In this example, we define $\Ccal := \{-1, \, 1\}$, and $N= 5\cdot 10^3$. The optimization is carried out on \verb+Matlab+ by using \verb+cvx+ toolbox.

The solution of Problem 1 for this example gives $\op{tr} \{(\Ical _F ^{\op{app}})^{-1}\} = 0.1108$ and $\det (\Ical _F ^{\op{app}}) = 1.8036 \cdot 10^{12}$. Figure \ref{fig: 4} presents the stationary probabilities of the optimal input signal for both cost functions. We can see that the stationary probabilities depend on the cost function $h$. However, we see that both cost functions assign higher stationary probabilities to the states $(-1,\, -1)$ and $(1,\, 1)$.

We can compare the performance of our approach with the method introduced in \cite{brighenti2009}. For this purpose, we generate an input sequence of length $N$ by running the Markov chain associated to the stationary distribution in Figure \ref{fig: 41}, and the 4-states Markov chain presented in \cite{brighenti2009}. To guarantee that the input is a realization of a stationary process, we discard the first $10^6$ outputs of the Markov chain. The results for the sampled information matrix are\footnote{Notice that our results are consistent with those reported in \cite{brighenti2009}, since the scaling factor $N$ is not considered here.} $\op{tr} \{(\Ical _F )^{-1}\} = 1.8233 \cdot 10^{-4}$ for the 4-states Markov chain presented in \cite{brighenti2009}, and $\op{tr} \{(\Ical _F )^{-1}\} = 1.6525 \cdot 10^{-4}$ for our method. Therefore, the approach in this paper gives better results for the example introduced in \cite{brighenti2009}.

To have an idea of the computation time required for this example, the optimization was solved in a laptop Dell Latitude E6430, equipped with Intel Core i7 $2.6$ [GHz] processor, and $8$ [Gb] of RAM memory. The time required from the computation of elementary cycles to the computation of stationary probabilities is $1.9$ seconds\footnote{A time bound for the computation of elementary cycles is given by $O(c_{\op{seq}}^{n_m} (c_{\op{seq}} +1)(c_e+1))$, where $c_e$ is the number of elementary cycles \cite[p. 77]{johnson1975}.}.
%Notice that we can improve the optimal value of both cost functions by increasing the length of the sequence $n_{\op{seq}}$. Indeed, if we solve this example for $n_{m}=3$ we obtain $\op{tr} \{(\Ical _F ^{\op{app}})^{-1}\} = 3.1367\cdot 10^{-2}$ and $\det (\Ical _F ^{\op{app}}) = 1.0036 \cdot 10^{14}$.
\fin

The numerical examples presented in this section show that the proposed method %proposed in this article
 is suitable for %solving
 input design %problems
 for systems with output-error-type structure, and when amplitude constraints on the input are required.
%%%%%%%%%%%%%%%%%%%%%%%%%%%%%%%%%%%%%%%%%%%%%%%%%%%%%%%%%%%%%%%%%%%%%%%%%%%%%%%%
\section{CONCLUSIONS}\label{sec: 5}
In this paper we have developed a new method to compute input signals for systems with arbitrary nonlinearities. The method is based on the optimization of a scalar cost function of the information matrix with respect to the probability density function of a stationary input. The optimal probability density function is used to compute the optimal input signal. %Since the problem is computationally intractable,
 An approach based on graph theory is used to derive a computationally efficient algorithm. %to solve the problem.
 This approach assumes that the input can adopt a finite set of values. An important feature of this method is that, by a suitable definition of the cost function,
 the optimization problem is convex even for nonlinear systems. Numerical examples show that this method is consistent with previous results in the literature, when we assume that the system has a particular structure. The method can also be used for input design with amplitude limitations. %, which is difficult to solve.

%Future work will be focused on input design for systems with more arbitrary structure, and design optimal transition probabilities for the Markov chain associated to the input sequence.
%%%%%%%%%%%%%%%%%%%%%%%%%%%%%%%%%%%%%%%%%%%%%%%%%%%%%%%%%%%%%%%%%%%%%%%%%%%%%%%%
\section{ACKNOWLEDGMENTS}\label{sec: 6}
The authors thank to Marco Forgione for the fruitful discussions during his time at the Automatic Control Lab in KTH.
%%%%%%%%%%%%%%%%%%%%%%%%%%%%%%%%%%%%%%%%%%%%%%%%%%%%%%%%%%%%%%%%%%%%%%%%%%%%%%%%
%%%%% References %%%%%
\bibliographystyle{IEEEbib}
\bibliography{library}%bib,bib2,bib3,bibliografia_PVP,bibliografia_EDU}

\begin{thebibliography}{10}

\bibitem{Cox1958}
D.R. Cox,
\newblock {\em Planning of experiments},
\newblock New York: Wiley, 1958.

\bibitem{goopay76}
G.C. Goodwin and R.L. Payne,
\newblock {\em {Dynamic System Identification: Experiment Design and Data
  Analysis}},
\newblock Academic Press, New York, 1977.

\bibitem{fedorov1972}
V.V. Fedorov,
\newblock {\em {Theory of optimal experiments}},
\newblock New York: Academic Press, 1972.

\bibitem{Whittle1973}
P.~Whittle,
\newblock ``Some general points in the theory of optimal exexperiment design,''
\newblock {\em Journal of {R}oyal {S}tatistical {S}ociety}, vol. 1, pp.
  123--130, 1973.

\bibitem{hildebrand2003}
R~Hildebrand and M.~Gevers,
\newblock ``Identification for control: {O}ptimal input design with respect to
  a worst-case {$\nu$}-gap cost function,''
\newblock {\em {SIAM} {J}ournal of {C}ontrol {O}ptimization}, vol. 41, no. 5,
  pp. 1586--1608, 2003.

\bibitem{gevers2005}
M.~Gevers,
\newblock ``{Identification for control: from the early achievements to the
  revival of experiment design},''
\newblock {\em European Journal of Control}, vol. 11, pp. 1--18, 2005.

\bibitem{jansson2005}
H.~Jansson and H.~Hjalmarsson,
\newblock ``Input design via {LMI}s admitting frequency-wise model
  specifications in confidence regions,''
\newblock {\em {IEEE} {T}ransactions on {A}utomatic {C}ontrol}, vol. 50, no.
  10, pp. 1534--1549, Oct. 2005.

\bibitem{lindqvist2000}
K.~Lindqvist and H.~Hjalmarsson,
\newblock ``Optimal input design using linear matrix inequalities,''
\newblock in {\em {IFAC} {S}ymposium on {S}ystem {I}dentification}, Santa
  Barbara, California, {USA}, July 2000.

\bibitem{brighenti2009}
C.~Brighenti, B.~Wahlberg, and C.R. Rojas,
\newblock ``Input design using {M}arkov chains for system identification,''
\newblock in {\em Joint 48th {C}onference on {D}ecision and {C}ontrol and 28th
  {C}hinese {C}onference}, Shangai, P.R. China, 2009, pp. 1557--1562.

\bibitem{Rojas2007}
C.R. Rojas, J.S. Welsh, G.C. Goodwin, and A.~Feuer,
\newblock ``{Robust optimal experiment design for system identification},''
\newblock {\em Automatica}, vol. 43, no. 6, pp. 993--1008, June 2007.

\bibitem{Suzuki2007}
H.~Suzuki and T.~Sugie,
\newblock ``On input design for system identification in time domain,''
\newblock in {\em Proceedings of the {E}uropean {C}ontrol {C}onference}, Kos,
  Greece, July 2007.

\bibitem{hjalmarsson2007}
H.~Hjalmarsson and J.~Mårtensson,
\newblock ``Optimal input design for identification of non-linear systems:
  {L}earning from the linear case,''
\newblock in {\em {A}merican {C}ontrol {C}onference}, New York, {U}nited
  {S}tates, 2007, pp. 1572--1576.

\bibitem{vincent2009}
T.L. Vincent, C.~Novara, K.~Hsu, and K.~Poola,
\newblock ``Input design for structured nonlinear system identification,''
\newblock in {\em 15th {IFAC} {S}ymposium on {S}ystem {I}dentification},
  Saint-{M}alo, {F}rance, 2009, pp. 174--179.

\bibitem{larsson2010}
C.~Larsson, H.~Hjalmarsson, and C.R. Rojas,
\newblock ``{On optimal input design for nonlinear FIR-type systems},''
\newblock in {\em 49th IEEE Conference on Decision and Control}, Atlanta, USA,
  2010, pp. 7220--7225.

\bibitem{gopaluni2011}
R.B. Gopaluni, T.B. Schön, and A.G. Wills,
\newblock ``Input design for nonlinear stochastic dynamic systems - {A}
  particle filter approach,''
\newblock in {\em 18th {IFAC} {W}orld {C}ongress}, {M}ilano, {I}taly, 2011.

\bibitem{zaman1983}
A.~Zaman,
\newblock ``{Stationarity on finite strings and shift register sequences},''
\newblock {\em The Annals of Probability}, vol. 11, no. 3, pp. 678--684, Aug.
  1983.

\bibitem{johnson1975}
D.B. Johnson,
\newblock ``{Finding all the elementary circuits of a directed graph},''
\newblock {\em SIAM Journal on Computing}, vol. 4, no. 1, pp. 77--84, March
  1975.

\bibitem{tarjan1972}
R.~Tarjan,
\newblock ``{Depth-First Search and Linear Graph Algorithms},''
\newblock {\em SIAM Journal on Computing}, vol. 1, no. 2, pp. 146--160, June
  1972.

\bibitem{rojas2011adaptive}
C.R. Rojas, H.~Hjalmarsson, L.~Gerencs{\'e}r, and J.~M{\aa}rtensson,
\newblock ``An adaptive method for consistent estimation of real-valued
  non-minimum phase zeros in stable {LTI} systems,''
\newblock {\em Automatica}, vol. 47, no. 7, pp. 1388--1398, 2011.

\bibitem{ljung99}
L.~Ljung,
\newblock {\em {System Identification. Theory for the User, 2nd ed.}},
\newblock Upper Saddle River, NJ: Prentice-Hall, 1999.

\bibitem{boyvan04}
S.~Boyd and L.~Vandenberghe,
\newblock {\em {Convex Optimization}},
\newblock Cambridge University Press, 2004.

\bibitem{hastings1970monte}
W.K. Hastings,
\newblock ``Monte {C}arlo sampling methods using {M}arkov chains and their
  applications,''
\newblock {\em Biometrika}, vol. 57, no. 1, pp. 97--109, 1970.

\bibitem{boyd2004}
S.~Boyd, P.~Diaconis, and L.~Xiao,
\newblock ``{Fastest mixing Markov chain on a graph},''
\newblock {\em SIAM Review}, vol. 46, no. 4, pp. 667--689, Oct. 2004.

\bibitem{boyd_cvxuser2013}
M.C. Grant and S.P. Boyd,
\newblock {\em The {CVX} users' guide},
\newblock {CVX} {R}esearch, {I}nc., 2nd. edition, January 2013.

\end{thebibliography}

%\begin{thebibliography}{99}
%
%\end{thebibliography}
\end{document}